\documentclass[a4paper,11pt]{amsart}
\usepackage{amssymb}
\usepackage[arrow,curve,matrix]{xy}

\title[Derived Equivalences of Finite Posets]
{On Derived Equivalences of Categories of Sheaves Over Finite
Posets}

\author{Sefi Ladkani}
\address{Einstein Institute of Mathematics, The Hebrew University of Jerusalem, Jerusalem 91904, Israel}
\email{sefil@math.huji.ac.il}

\newcommand{\one}{\mathbf{1}}
\newcommand{\ii}{i_{*}i^{-1}}
\newcommand{\jj}{j_{!}j^{-1}}
\newcommand{\wtP}{\widetilde{P}}
\newcommand{\wtI}{\widetilde{I}}

\newcommand{\cA}{\mathcal{A}}
\newcommand{\cD}{\mathcal{D}}
\newcommand{\cE}{\mathcal{E}}
\newcommand{\cF}{\mathcal{F}}
\newcommand{\cG}{\mathcal{G}}
\newcommand{\cI}{\mathcal{I}}
\newcommand{\cK}{\mathcal{K}}

\newcommand{\cT}{\mathcal{T}}
\newcommand{\cX}{\mathcal{X}}
\newcommand{\bZ}{\mathbb{Z}}
\newcommand{\bQ}{\mathbb{Q}}
\newcommand{\gL}{\Lambda}

\newcommand{\dA}{\cD^b(A)}
\newcommand{\dX}{\cD^b(X)}
\newcommand{\dY}{\cD^b(Y)}
\newcommand{\dU}{\cD^b(U)}

\newcommand{\fX}{\mathfrak{X}}

\DeclareMathOperator{\GL}{GL}
\DeclareMathOperator{\supp}{supp}
\DeclareMathOperator{\gldim}{gl.\!dim}
\DeclareMathOperator{\Hom}{Hom}
\DeclareMathOperator{\Ext}{Ext}
\DeclareMathOperator{\End}{End}
\DeclareMathOperator{\hh}{H}
\DeclareMathOperator{\HH}{HH}
\DeclareMathOperator{\Id}{Id}

\newcommand{\bform}[2]{\left\langle {#1}, {#2} \right\rangle}

\theoremstyle{plain}
\newtheorem{theorem}{Theorem}[section]
\newtheorem{lemma}[theorem]{Lemma}
\newtheorem{prop}[theorem]{Proposition}
\newtheorem{cor}[theorem]{Corollary}

\theoremstyle{definition}
\newtheorem{defn}[theorem]{Definition}

\newtheorem{example}[theorem]{Example}

\numberwithin{equation}{section}

\begin{document}

\begin{abstract}
A finite poset $X$ carries a natural structure of a topological
space. Fix a field $k$, and denote by $\dX$ the bounded derived
category of sheaves of finite dimensional $k$-vector spaces over
$X$. Two posets $X$ and $Y$ are said to be \emph{derived equivalent}
if $\dX$ and $\dY$ are equivalent as triangulated categories.

We give explicit combinatorial properties of $X$ which are invariant
under derived equivalence, among them are the number of points, the
$\bZ$-congruency class of the incidence matrix, and the Betti
numbers. We also show that taking opposites and products preserves
derived equivalence.

For any closed subset $Y \subseteq X$, we construct a strongly
exceptional collection in $\dX$ and use it to show an equivalence
$\dX \simeq \dA$ for a finite dimensional algebra $A$ (depending on
$Y$). We give conditions on $X$ and $Y$ under which $A$ becomes an
incidence algebra of a poset.

We deduce that a lexicographic sum of a collection of posets along a
bipartite graph $S$ is derived equivalent to the lexicographic sum
of the same collection along the opposite $S^{op}$.

This construction produces many new derived equivalences of posets
and generalizes other well known ones.

As a corollary we show that the derived equivalence class of an
ordinal sum of two posets does not depend on the order of summands.
We give an example that this is not true for three summands.
\end{abstract}

\maketitle

\section{Introduction}

Over the last years a growing interest in the understanding of
derived categories of coherent sheaves over algebraic varieties, and
in particular, the question when two varieties have equivalent
derived categories of sheaves, has emerged~\cite{BondalOrlov02}.

We investigate a similar question for partially ordered sets
(\emph{posets}). A poset $X$ carries a natural structure of a
topological space, therefore we can consider the category of sheaves
over $X$ with values in an abelian category $\cA$.

We focus on the case where $\cA$ is the category of finite
dimensional vector spaces over a field $k$, which allows us to
identify the category of sheaves with a category of modules over the
incidence algebra of $X$ over $k$, so that tools from the theory of
derived equivalence of algebras can be used. However, there is no
known algorithm which decides, given two posets, whether their
derived categories of sheaves of finite dimensional $k$-vector
spaces are equivalent.

In Section~2, we present in a specific way, appropriate for dealing
with posets, the basic notions from sheaf theory that will be used
throughout the paper. In Section~3 we discuss combinatorial
invariants of derived equivalence, whereas in Section~4 we
construct, for any poset $X$ admitting a special structure, new
poset derived equivalent to $X$. This construction is based on the
notion of strongly exceptional sequences in triangulated categories
and partially generalizes the known constructions
of~\cite{APR79,BGP73}.

\subsection*{Acknowledgement}
I am grateful to T.~Holm, B.~Keller and my advisor D.~Kazhdan for
helpful discussions related to this paper.

\section{Preliminaries}
\subsection{Finite posets and $T_0$ spaces}
Throughout this note, the term \emph{poset} will mean a finite
partially ordered set. Any poset $(X, \leq)$ carries a structure of
a topological space by defining the closed sets to be the subsets $Y
\subseteq X$ such that if $y \in Y$ and $y' \leq y$ then $y' \in Y$.

For each $x \in X$, denote by $\{x\}^{-}$ the closure of $\{x\}$ and
by $U_x$ the minimal open subset of $X$ containing $x$, which equals
the intersection of the open sets containing $x$. Then $\{x\}^{-} =
\left\{ x' \in X \,:\, x' \leq x \right\}$, $U_x = \left\{ x' \in X
\,:\, x' \geq x \right\}$ and
\[
x \leq x' \Longleftrightarrow \{x\}^- \subseteq \{x'\}^-
\Longleftrightarrow U_{x'} \subseteq U_x
\]

If $x, y$ are two distinct points in $X$, then one of the open sets
$U_x, U_y$ does not contain both points, thus $X$ satisfies the
$T_0$ separation property.

Conversely, given a finite $T_0$ topological space $X$, let $U_x$ be
the intersection of all open sets in $X$ containing $x \in X$.
Define a partial order $\leq$ on $X$ by $x \leq x'$ if $U_{x'}
\subseteq U_x$.

This leads to an identification of posets with finite $T_0$
topological spaces. Such spaces have been studied in the past
\cite{McCord66,Stong66}, where it turned out that their homotopy and
homology properties are more interesting than might seem at first
glance. For example, if $\cK$ is any finite simplicial complex and
$X$ is the $T_0$ space induced by the partial order on the simplices
of $\cK$, then there exists a weak homotopy equivalence $|\cK| \to
X$ \cite{McCord66}.

\subsection{Sheaves and diagrams}
Given a poset $X$, its \emph{Hasse diagram} is a directed graph
defined as follows. Its vertices are the elements of $X$ and its
directed edges $x \to y$ are the pairs $x < y$ in $X$ such that
there is no $z \in X$ with $x < z < y$. The anti-symmetry condition
on $\leq$ implies that this graph has no directed cycles.

Let $X$ be a poset and $\cA$ be an abelian category. Using the
topology on $X$, we can consider the category of \emph{sheaves over
$X$ with values in $\cA$}, denoted by $Sh_X \cA$ or sometimes
$\cA^X$.

We note that sheaves over posets were used in~\cite{DGM00} for the
computation of cohomologies of real subspace arrangements.

Let $\cF$ be a sheaf on $X$. If $x \in X$, let $\cF(x)$ be the stalk
of $\cF$ over $x$, which equals $\cF(U_x)$. The restriction maps
$\cF(x) = \cF(U_x) \to \cF(U_{x'}) =\cF(x')$ for $x' > x$ give rise
to a commutative diagram over the Hasse diagram of $X$. Conversely,
such a diagram $\{F_x\}$  defines a sheaf $\cF$ by setting the
sections as the inverse limits $\cF(U) = \lim_{x \in U} F_x$.
Indeed, it is enough to verify the sheaf condition for the sets
$U_x$, which follows from the observation that for any cover $U_x =
\bigcup_i U_{z_i}$, one of the $z_i$ equals $x$.

Thus we may identify $Sh_X \cA$ with the category of commutative
diagrams over the Hasse diagram of $X$ and interchange the terms
sheaf and diagram as appropriate. The latter category can be viewed
as the category of functors $X \to \cA$ where we consider $X$ as a
category whose objects are the points $x \in X$, with unique
morphisms $x \to x'$ for $x \leq x'$. Under this identification, the
global sections functor $\Gamma(X; -) : \cA^X \to \cA$ defined as
$\Gamma(X; \cF) = \cF(X)$, coincides with the (inverse) limit
functor $\lim_X : \cA^X \to \cA$.

\subsection{Functors associated with a map $f:X \to Y$}
A map $f:X \to Y$ between two finite posets is continuous if and
only if it is \emph{order preserving}, that is, $f(x) \leq f(x')$
for any $x \leq x'$ in $X$ \cite[Prop.~7]{Stong66}.

A continuous map $f:X \to Y$ gives rise to the functors $f_{*},
f_{!}: Sh_X \cA \to Sh_Y \cA$ and $f^{-1}: Sh_Y \cA \to Sh_X \cA$
defined, in terms of diagrams, by
\begin{align*}
(f^{-1} \cG)(x) &= \cG(f(x)) \\
(f_{*} \cF)(y) &= \lim_{\xleftarrow{}} \{\cF(x) \,:\, f(x) \geq y \}
\\
(f_{!} \cF)(y) &= \lim_{\xrightarrow{}} \{\cF(x) \,:\, f(x) \leq y\}
\end{align*}
where $x \in X$, $y \in Y$ and $\cF \in Sh_X \cA$, $\cG \in Sh_Y
\cA$. Viewing $X$, $Y$ as categories and $\cF \in Sh_X \cA$ as a
functor $\cF : X \to \cA$, the sheaves $f_{*} \cF$ and $f_{!} \cF$
are the right and left Kan extensions of $\cF$ along $f: X \to Y$.

The functors $f^{-1}, f_{*}$ coincide with the usual ones from sheaf
theory. We have the following adjunctions:
\begin{align}
\label{e:adj}
\Hom_{Sh_X \cA}(f^{-1}\cG, \cF) &\simeq \Hom_{Sh_Y \cA}(\cG, f_{*}\cF) \\
\notag \Hom_{Sh_X \cA}(\cF, f^{-1}\cG) &\simeq \Hom_{Sh_Y
\cA}(f_{!}\cF, \cG)
\end{align}
so that $f_{*}$ is left exact and $f_{!}$ is right exact. $f^{-1}$
is exact, as can be seen from its action on the stalks.

If $Y$ is a closed subset of $X$, we have a \emph{closed embedding}
$i:Y \to X$. In this case, $i_{*}$ is exact. This is because $i_{*}$
takes a diagram on $Y$ and extends it to $X$ by filling the vertices
of $X \setminus Y$ with zeros. Similarly, for an open embedding $j:U
\to X$, $j_{!}$ is exact, as it extends by zeros diagrams on $U$.
Now let $Y \subseteq X$ be closed and $U = X \setminus Y$ its
complement. The adjunction morphisms $\jj \cF \to \cF$ and $\cF \to
\ii \cF$ for the embeddings $i: Y \to X$ and $j: U \to X$ induce a
short exact sequence
\begin{equation}
\label{e:jjii} 0 \to j_{!}j^{-1} \cF \to \cF \to i_{*}i^{-1} \cF \to
0
\end{equation}
for any sheaf $\cF$ on $X$, as can be verified at the stalks.

\subsection{Simples, projectives and injectives}
When $f:X \to \bullet$ is the mapping to a point, $f_{*} = \Gamma(X;
-)$ and $f^{-1}(M)$ for an object $M$ of $\cA$ gives the
\emph{constant sheaf} on $X$ with value $M$.

Let $x \in X$ and consider the map $i_x : \bullet \to X$ whose image
is $\{x\}$. Then $i_x^{-1}(\cF) = \cF(x)$ is the stalk at $x$ and
for an object $M$ of $\cA$ we have
\begin{align*}
({i_x}_* M)(y) = \begin{cases} M & \text{if $y \leq x$} \\
0 & \text{otherwise}
\end{cases}
& &
({i_x}_! M)(y) = \begin{cases} M & \text{if $y \geq x$} \\
0 & \text{otherwise}
\end{cases}
\end{align*}
with identity arrows between the $M$-s. The
adjunctions~\eqref{e:adj} take the form:
\begin{align}
\label{e:adjx}
&\Hom_{Sh_X \cA}(\cF, {i_x}_{*} M) \simeq \Hom_{\cA}(\cF(x), M) \\
\notag &\Hom_{Sh_X \cA}({i_x}_{!} M, \cF) \simeq \Hom_{\cA}(M,
\cF(x))
\end{align}
and we deduce the following lemma:
\begin{lemma}
If $I$ is injective in $\cA$, ${i_x}_* I$ is injective in $Sh_X
\cA$. If $P$ is projective in $\cA$, ${i_x}_! P$ is projective in
$Sh_X \cA$.
\end{lemma}
\begin{cor}
\label{c:injproj} If $\cA$ has enough injectives (projectives), so
does $Sh_X \cA$.
\end{cor}
\begin{proof}
The identity maps $\cF(x) \xrightarrow{=} \cF(x)$ induce,
via~\eqref{e:adjx}, an injection $\cF \hookrightarrow \oplus_{x \in
X} {i_x}_* \cF(x)$ and surjection $\oplus_{x \in X} {i_x}_! \cF(x)
\twoheadrightarrow \cF$. Now replace each $\cF(x)$ by an injective
(or projective) cover.
\end{proof}

For a sheaf $\cF$, let $\supp \cF = \{ x \in X \,:\, \cF(x) \neq 0
\}$ be its \emph{support}. We call $\cF$ a \emph{stalk sheaf} if its
support is a point. For any object $M$ of $\cA$ and $x \in X$ there
exists a stalk sheaf $M_x$ whose stalk at $x$ equals $M$. Moreover
$M_x$ is simple in $Sh_X \cA$ if and only if $M$ is simple in $\cA$.

The following lemma is proved by induction on the number of elements
$|X|$, using~\eqref{e:jjii} and the fact that the partial order on
$X$ can be extended to a linear order, i.e. one can write the
elements of $X$ in a sequence $x_1,x_2,\dots,x_n$ such that for any
$1 \leq i,j \leq n$, $x_i < x_j$ implies that $i < j$.

\begin{lemma}
\label{l:filt} Any sheaf $\cF$ on $X$ admits a finite filtration
whose quotients are stalk sheaves.
\end{lemma}

Denote by $\gldim \cA$ the \emph{global dimension} of an abelian
category $\cA$. This is the maximal integer $n$ for which there
exist objects $M, M'$ of $\cA$ with $\Ext^n(M, M') \neq 0$ (and
$\infty$ if there is no such maximal $n$). Recall that an abelian
category is a \emph{finite length} category if every object is of
finite length. From Lemma~\ref{l:filt}, we have:
\begin{cor}
If $\cA$ is a finite length category, so is $Sh_X \cA$.
\end{cor}

\begin{defn}
A strictly increasing sequence $x_0 < x_1 < \dots < x_n$ in $X$ is
called a \emph{chain of length $n$}. The \emph{dimension} of $X$,
denoted $\dim X$, is the maximal length of a chain in $X$.
\end{defn}

\begin{prop}[\cite{Mitchell68}]
\label{p:gldim} $\gldim Sh_X \cA \leq \gldim \cA + \dim X$.
\end{prop}

The difference $\gldim Sh_X \cA - \gldim \cA$ obviously depends on
$X$, but it may well depend also on $\cA$, see the examples in
\cite{IgusaZacharia90,Spears72}.


\subsection{Sheaves of finite-dimensional vector spaces}
Fix a field $k$ and consider the category $\cA$ of finite
dimensional vector spaces over $k$. Denote by $Sh_X$ the category
$Sh_X \cA$ and by $\Hom_X(-,-)$ the morphism spaces
$\Hom_{Sh_X}(-,-)$ (We omit the reference to $k$ to emphasize that
it is to be fixed throughout).

The \emph{incidence algebra} of $X$ over $k$, denoted $kX$, is the
algebra spanned by $e_{xy}$ for the pairs $x \leq y$ in $X$, with
multiplication defined by $e_{xy}e_{zw} = \delta_{yz} e_{xw}$.

\begin{lemma}
The category $Sh_X$ is equivalent to the category of finite
dimensional right modules over the incidence algebra $kX$.
\end{lemma}
\begin{proof}
The proof is similar to the corresponding fact about representations
of a quiver and right modules over its path algebra. Namely, for a
sheaf $\cF$, consider $M = \oplus_{x \in X} \cF(x)$ and let $\iota_x
: \cF(x) \to M$, $\pi_x : M \to \cF(x)$ be the natural maps. Equip
$M$ with a structure of a right $kX$-module by letting the basis
elements $e_{xx'}$ for $x \leq x'$ act from the right as the
composition $M \xrightarrow{\pi_x} \cF(x) \to \cF(x')
\xrightarrow{\iota_{x'}} M$. Conversely, given a finite dimensional
right module $M$ over $kX$, set $\cF(x) = M e_{xx}$ and define the
maps $\cF(x) \to \cF(x')$ using the right multiplication by
$e_{xx'}$.
\end{proof}

The one dimensional space $k$ is both simple, projective and
injective in the category of $k$-vector spaces. Applying the results
of the previous subsection, we get, for any $x \in X$, sheaves $S_x,
P_x, I_x$ which are simple, projective and injective, respectively.
Explicitly,
\[
S_x(y) = \begin{cases} k & y=x \\ 0 & \text{otherwise} \end{cases}
,\, P_x(y) = \begin{cases} k & y \geq x \\ 0 & \text{otherwise}
\end{cases} ,\, I_x(y) = \begin{cases} k & y \leq x \\ 0 &
\text{otherwise}
\end{cases}
\]

By~\eqref{e:adjx}, for any sheaf $\cF$, $\Hom_X(P_x, \cF) = \cF(x)$
and $\Hom_X(\cF, I_x) = \cF(x)^{\vee}$ (the dual space). Since the
sets $U_x, \{x\}^{-}$ are connected, the sheaves $P_x, I_x$ are
indecomposable. The sheaves $S_x, P_x, I_x$ form a complete set of
representatives of the isomorphism classes of simples,
indecomposable projectives and indecomposable injectives
(respectively) in $kX$.

By Corollary~\ref{c:injproj}, $Sh_X$ has enough projectives and
injectives (note that this can also be deduced by its identification
with the category of finite dimensional modules over a finite
dimensional algebra). It has finite global dimension, since by
Proposition~\ref{p:gldim}, $\gldim Sh_X \leq \dim X$.

\begin{prop}
$Sh_X$ and $Sh_Y$ are equivalent if and only if $X$ and $Y$ are
isomorphic (as posets).
\end{prop}
\begin{proof}
Since the isomorphism classes of simple objects in $Sh_X$ are in
one-to-one correspondence with the elements $x \in X$, and for two
such simples $S_x, S_y$, $\dim_k \Ext^1 (S_x, S_y)$ equals $1$ if
there is a directed edge $x \to y$ in the Hasse diagram of $X$ and
$0$ otherwise, we see that the Hasse diagram of $X$, hence $X$, can
be recovered (up to isomorphism) from the category $Sh_X$.
\end{proof}

\subsection{The derived category of sheaves over a poset}
For a poset $X$, denote by $\dX$ the bounded derived category of
$Sh_X$.

If $\cE$ is a set of objects of a triangulated category $\cT$, we
denote by $\langle \cE \rangle$ the triangulated subcategory of
$\cT$ \emph{generated by $\cE$}, that is, the minimal triangulated
subcategory containing $\cE$. We say that $\cE$ \emph{generates}
$\cT$ if $\langle \cE \rangle = \cT$.

Since $Sh_X$ is of finite global dimension with enough projectives
and injectives, $\dX$ can be identified with the homotopy category
of bounded complexes of projectives (or bounded complexes of
injectives). Hence the collections $\{P_x\}_{x \in X}$ and
$\{I_x\}_{x \in X}$ generate $\dX$.

\begin{lemma}
\label{l:PxPy} Let $x, y \in X$ and $i \in \bZ$. Then
\[
\Hom_{\dX}(P_x, P_y[i]) = \Hom_{\dX}(I_x, I_y[i]) = \begin{cases} k
& \text{$y \leq x$ and $i=0$} \\ 0 & \text{otherwise}
\end{cases}
\]
\end{lemma}
\begin{proof}
Since $P_x$ is projective, $\Hom_{\dX}(P_x, \cF[i]) = 0$ for any
sheaf $\cF$ and $i \neq 0$. If $x, y \in X$, then
\[
\Hom_{\dX}(P_x, P_y) = \Hom_X(P_x, P_y) = P_y (x) = \begin{cases} k
& \text{if $x \geq y$} \\ 0 & \text{otherwise} \end{cases}
\]
The proof for $\{I_x\}_{x \in X}$ is similar.
\end{proof}

For a continuous map $f: X \to Y$, denote by $Rf_{*}, Lf_{!},
f^{-1}$ the derived functors of $f_*, f_!, f^{-1}$. The
adjunctions~\eqref{e:adj} imply that
\begin{align}
\label{e:adjder}
\Hom_{\dX}(f^{-1}\cG, \cF) &\simeq \Hom_{\dY}(\cG, Rf_{*}\cF) \\
\notag \Hom_{\dX}(\cF, f^{-1}\cG) &\simeq \Hom_{\dY}(Lf_{!}\cF, \cG)
\end{align}
for $\cF \in \dX$, $\cG \in \dY$.

\begin{defn}
We say that two posets $X$ and $Y$ are \emph{derived equivalent},
denoted $X \sim Y$, if the categories $\dX$ and $\dY$ are equivalent
as triangulated categories.
\end{defn}

\section{Combinatorial invariants of derived equivalence}
We give a list of combinatorial properties of posets which are
preserved under derived equivalence. Most of the properties are
deduced from known invariants of derived categories. For the
convenience of the reader, we review the relevant definitions.

\subsection{The number of points and $K$-groups}
Recall that for an abelian category $\cA$, the \emph{Grothendieck
group} $K_0(\cA)$ is the quotient of the free abelian group
generated by the isomorphism classes $[X]$ of objects $X$ of $\cA$
divided by the subgroup generated by the expressions $[X] - [Y] +
[Z]$ for all the short exact sequences $0 \to X \to Y \to Z \to 0$
in $\cA$.

Similarly, for a triangulated category $\cT$, the group $K_0(\cT)$
is the quotient of the free abelian group on the isomorphism classes
of objects of $\cT$ divided by its subgroup generated by $[X] - [Y]
+ [Z]$ for all the triangles $X \to Y \to Z \to X[1]$ in $\cT$
(where $[1]$ denotes the shift). The natural inclusion $\cA \to
\cD^b(\cA)$ induces an isomorphism $K_0(\cA) \cong K_0(\cD^b(\cA))$.

Let $X$ be a poset and denote by $|X|$ the number of points of $X$.
Denote by $K_0(X)$ the group $K_0(\dX)$.

\begin{prop}
$K_0(X)$ is free abelian of rank $|X|$.
\end{prop}
\begin{proof}
The set $\{S_x\}_{x \in X}$ forms a complete set of representatives
of the isomorphism classes of simple finite dimensional
$kX$-modules, hence it is a $\bZ$-basis of $K_0(X)$ (alternatively
one could use the filtration of Lemma~\ref{l:filt}).
\end{proof}

\begin{cor}
If $X \sim Y$ then $|X| = |Y|$.
\end{cor}

It is known \cite{DuggerShipley04} that rings with equivalent
derived categories have the same $K$-theory. However, higher
$K$-groups do not lead to refined invariants of the number of
points.

\begin{prop}
$K_i(Sh_X) \simeq K_i(Sh_\bullet)^{|X|}$ for $i \geq 0$.
\end{prop}
\begin{proof}
$Sh_X$ is a finite length category and by~\cite[Corollary~1,
p.~104]{Quillen73},
\[
K_i(Sh_X) \simeq \bigoplus_{x \in X} K_i(\End_X(S_x))
\]
Clearly, $k = \End_X(S_x)$.
\end{proof}

\subsection{Connected components}
For two additive categories $\cT_1$, $\cT_2$, consider the category
$\cT = \cT_1 \times \cT_2$ whose objects are pairs $(M_1, M_2)$ and
the morphisms are defined by
\[
\Hom_{\cT}((M_1, M_2), (N_1, N_2)) = \Hom_{\cT_1}(M_1, N_1) \times
\Hom_{\cT_2}(M_2, N_2)
\]
$\cT_1$, $\cT_2$ are embedded in $\cT$ via the fully faithful
functors $M_1 \mapsto (M_1, 0)$ and $M_2 \mapsto (0, M_2)$. Denoting
the images again by $\cT_1, \cT_2$, we have that $\Hom_{\cT}(\cT_1,
\cT_2) = 0$. In addition, the indecomposables in $\cT$ are of the
form $(M_1, 0)$ or $(0, M_2)$ for indecomposables $M_1 \in \cT_1$,
$M_2 \in \cT_2$.

An additive category $\cT$ is \emph{connected} if for any
equivalence $\cT \simeq \cT_1 \times \cT_2$, one of $\cT_1$, $\cT_2$
is zero.

\begin{defn}
A poset $X$ is \emph{connected} if it is connected as a topological
space. This is equivalent to the following
condition~\cite[Prop.~5]{Stong66}:

For any $x, y \in X$ there exists a sequence $x=x_0, x_1, \dots, x_n
= y$ in $X$ such that for all $0 \leq i < n$, either $x_i \leq
x_{i+1}$ or $x_i \geq x_{i+1}$.
\end{defn}

\begin{lemma}
If $X$ is connected then the category $\dX$ is connected.
\end{lemma}
\begin{proof}
Let $\dX \simeq \cT_1 \times \cT_2$ be an equivalence and consider
the indecomposable projectives $\{P_x\}_{x \in X}$. Since each $P_x$
is indecomposable, its image lies in $\cT_1$ or in $\cT_2$, and we
get a partition $X = X_1 \sqcup X_2$.

Assume that $X_1$ is not empty. Since $\Hom(P_x, P_y) \neq 0$ for
all $y \leq x$ and $\Hom(\cT_1, \cT_2) = 0$, $X_1$ must be both open
and closed in $X$, and by connectivity, $X_1 = X$. Moreover,
$\{P_x\}_{x \in X}$ generates $\dX$ as a triangulated category,
hence $\dX \simeq \cT_1$ and $\cT_2 = 0$.
\end{proof}

\begin{prop}
Let $X$ and $Y$ be two posets with decompositions
\begin{align*}
X = X_1 \sqcup X_2 \sqcup \dots \sqcup X_t & & Y = Y_1 \sqcup Y_2
\sqcup \dots \sqcup Y_s
\end{align*}
into connected components. If $X \sim Y$ then $s = t$ and there
exists a permutation $\pi$ on $\{1,\dots,s\}$ such that $X_i \sim
Y_{\pi(i)}$ for all $1 \leq i \leq s$.
\end{prop}
\begin{proof}
There exists a pair of equivalences
\[
\xymatrix{{\cD^b(X_1) \times \dots \times \cD^b(X_t) = \dX}
\ar@<1ex>[r]^F & {\dY = \cD^b(Y_1) \times \dots \times \cD^b(Y_s)}
\ar@<1ex>[l]^G}
\]

If $x \in X$, the image $F(P_x)$ is indecomposable in $\dY$, hence
lands in one of the $\cD^b(Y_j)$, and we get a function $f: X \to
\{1,\dots,s\}$. For any $x' \leq x$, $\Hom_X(P_x, P_{x'}) \neq 0$,
therefore $f$ is constant on the connected components $X_i$ and
induces a map $\pi_F : \{1,\dots,t\} \to \{1,\dots,s\}$ via
$\pi_F(i) = f(x)$ for $x \in X_i$. Moreover, since $\{P_x\}_{x \in
X_i}$ generates $\cD^b(X_i)$ as a triangulated category, $F$
restricts to functors $\cD^b(X_i) \to \cD^b(Y_{\pi_F(i)})$, $1 \leq
i \leq t$.

Similarly for $G$, we obtain a map $\pi_G : \{1, \dots, s\} \to \{1,
\dots, t\}$ and functors $\cD^b(Y_j) \to \cD^b(Y_{\pi_G(j)})$ which
are restrictions of $G$.

For any $1 \leq i \leq t$, the image of $\cD^b(X_i)$ under $GF$ lies
in $\cD^b(X_{\pi_G\pi_F(i)})$. Since $GF$ is isomorphic to the
identity functor but on the other hand there are no nonzero maps
between $\cD^b(X_i)$ and $\cD^b(X_{i'})$ for $i \neq i'$ (as we
think of $X_i$ as \emph{subsets} of $X$, not just as abstract
sets!), we get that $\pi_G\pi_F(i)=i$ so that $\pi_G\pi_F$ is
identity. Similarly, $\pi_F\pi_G$ is identity.

We deduce that $s=t$, $\pi_F$ and $\pi_G$ are permutations, and the
restrictions of $F$ induce equivalences $\cD^b(X_i) \simeq
\cD^b(Y_{\pi_F(i)})$.
\end{proof}

One can also deduce that the number of connected components is a
derived invariant by considering the center $Z(kX)$ of the incidence
algebra $kX$ using the fact that derived equivalent algebras have
isomorphic centers~\cite{Rickard89}.

\begin{lemma}
$Z(kX) \cong k \times k \times \dots  \times k$ where the number of
factors equals the number of connected components of $X$.
\end{lemma}
\begin{proof}
Let $c = \sum_{x \leq y} c_{xy} e_{xy} \in Z(kX)$. Comparison of the
coefficients of $e_{xx} c$ and $c e_{xx}$ gives $c_{xy}=0$ for $x
\neq y$, thus $c = \sum_{x} c_x e_{xx}$.

If $x \leq y$ then $c_x e_{xy} = c e_{xy} = e_{xy} c = c_y e_{xy}$,
hence $c_x = c_y$ if $x, y$ are in the same connected component.
\end{proof}

\subsection{The Euler form and M\"obius function}
Let $X$ be a poset. Since $Sh_X$ has finite global dimension, the
expression
\[
\bform{K}{L}_X = \sum_{i \in \bZ} (-1)^i \dim_k \Hom_{\dX}(K, L[i])
\]
is well-defined for $K, L \in \dX$ and induces a $\bZ$-bilinear form
on $K_0(X)$, known as the \emph{Euler form}.

Recall that the \emph{incidence matrix} of $X$, denoted $\one_X$, is
the $X \times X$ matrix defined by
\[
(\one_X)_{xy} = \begin{cases} 1 & x \leq y \\ 0 & \text{otherwise}
\end{cases}
\]
By extending the partial order on $X$ to a linear order, we can
always arrange the elements of $X$ such that the incidence matrix is
upper triangular with ones on the diagonal. In particular, $\one_X$
is invertible over $\bZ$.

\begin{defn}
The \emph{M\"obius function} $\mu_X : X \times X \to \bZ$ is defined
by $\mu_X(x, y) = (\one_X^{-1})_{xy}$.
\end{defn}

The following is an immediate consequence of the definition.
\begin{lemma}[M\"obius inversion formula]
\label{l:Mobius} Let $f : X \to \bZ$. Define $g: X \to \bZ$ by $g(x)
= \sum_{y \geq x} f(y)$. Then $f(x) = \sum_{y \geq x} \mu_X(x,y)
g(y)$.
\end{lemma}

The M\"obius inversion formula can be used to compute the matrix of
the Euler form with respect to the basis of simple objects.

\begin{lemma}
\label{l:PxSy} $\bform{[P_x]}{[S_y]}_X = \delta_{xy}$ for all $x, y
\in X$.
\end{lemma}
\begin{proof}
Since $P_x$ is projective, $\Hom_{\dX}(P_x, \cF[i]) = 0$ for any
sheaf $\cF$ and $i \neq 0$. Now by~\eqref{e:adjx}, $\Hom_{\dX}(P_x,
S_y) = \Hom_X(P_x, S_y) = S_y(x)$.
\end{proof}

\begin{prop}
\label{p:Euler} Let $x, y \in X$. Then $\bform{[S_x]}{[S_y]}_X =
\mu_X(x, y)$.
\end{prop}
\begin{proof}
Fix $y$ and define $f:X \to \bZ$ by $f(x) = \bform{[S_x]}{[S_y]}_X$.
Since $[P_x] = \sum_{x' \geq x} [S_{x'}]$, Lemmas~\ref{l:Mobius}
and~\ref{l:PxSy} imply that
\[
f(x) = \sum_{x' \geq x} \mu_X(x,x') \bform{[P_{x'}]}{[S_y]}_X =
\mu_X(x,y)
\]
\end{proof}

\begin{defn}
Let $R$ be a commutative ring. Two matrices $M_1, M_2 \in \GL_n(R)$
are \emph{congruent} over $R$ if there exists a matrix $P \in
\GL_n(R)$ such that $M_2 = P M_1 P^t$.
\end{defn}

Note that if $M_1, M_2$ is a pair of congruent matrices, so are
$M_1^t, M_2^t$ and $M_1^{-1}, M_2^{-1}$. Denote by $M^{-t}$ the
inverse of the transpose of $M$.

\begin{cor}
If $X \sim Y$ then $\one_X$, $\one_Y$ are congruent over $\bZ$.
\end{cor}
\begin{proof}
An equivalence $F : \dX \to \dY$ induces an isomorphism $[F] :
K_0(X) \to K_0(Y)$ which preserves the Euler form. By
Proposition~\ref{p:Euler}, the matrix of the Euler form of $\dX$
over the basis of simples is $\one_X^{-1}$, hence $[F]^t \one_Y^{-1}
[F] = \one_X^{-1}$.
\end{proof}

In practice, testing for congruence over $\bZ$ is not an easy task.
However, the following necessary condition is often very useful in
ruling out congruence.

\begin{lemma}
Let $M_1, M_2 \in \GL_n(R)$ be congruent. Then the matrices $M_1
M_1^{-t}$, $M_2 M_2^{-t}$ are conjugate in the group $\GL_n(R)$.
\end{lemma}
\begin{proof}
If $M_2 = P M_1 P^t$ for some $P \in \GL_n(R)$, then
\[
M_2 M_2^{-t} = (P M_1 P^t) (P^{-t} M_1^{-t} P^{-1}) = P M_1 M_1^{-t}
P^{-1}
\]
\end{proof}

\begin{cor}
\label{c:simQFp} If $X \sim Y$ then $\one_X \one_X^{-t}$ and $\one_Y
\one_Y^{-t}$ are similar over $\bZ$. In particular, they are similar
over $\bQ$ and modulo all primes $p$.
\end{cor}

Note that $\one_X \one_X^{-t}$ is (up to sign) the \emph{Coxeter
matrix} of the algebra $kX$. It is the image in $K_0(X)$ of the
Serre functor on $\dX$.

\subsection{Betti numbers and Euler characteristic}
The Hochschild cohomology is a known derived invariant of an
algebra~\cite{Happel88,Rickard91}. For posets, one can compute the
Hochschild cohomology as the simplicial cohomology of an appropriate
simplicial complex~\cite{Cibils89,GerstenhaberSchack83}. Thus the
simplicial cohomology is a derived invariant, which we relate to the
cohomology of the constant sheaf.

For the convenience of the reader, we review the notions of sheaf
cohomology, simplicial cohomology and Hochschild cohomology. As
before, we keep the field $k$ fixed.

\subsubsection{Sheaf cohomology}
Recall that the $i$-th cohomology of a sheaf $\cF \in Sh_X$, denoted
$\hh^i(X;\cF)$, is the value of the $i$-th right derived functor of
the global sections functor $\Gamma(X;-) : Sh_X \to Sh_{\bullet}$.
Observe that $\Gamma(X;\cF) = \Hom_X(k_X, \cF)$ where $k_X$ is the
\emph{constant sheaf} on $X$, i.e. $k_X(x) = k$ for all $x \in X$
with all morphisms being the identity of $k$. It follows that
$\hh^i(X;\cF) = \Ext^i_X(k_X, \cF)$. Specializing this for the
particularly interesting cohomologies of the constant sheaf, we get
that $\hh^i(X; k_X) = \Ext^i_X(k_X, k_X)$.

\subsubsection{Simplicial cohomology}
Let $X$ be a poset, $p \geq 0$. A \emph{$p$-dimensional simplex} in
$X$ is a chain of length $p$. Since subsets of chains are again
chains, the set of all simplices in $X$ forms a simplicial complex
$\cK(X)$~\cite{McCord66}, known as the \emph{order complex} of $X$.
The $i$-th simplicial cohomology of $X$ is defined as the $i$-th
simplicial cohomology of $\cK(X)$, and we denote it by $\hh^i(X)$.
The number $\beta^i(X) = \dim_k \hh^i(X)$ is the $i$-th \emph{Betti
number} of $X$.

The simplicial cohomology of $X$ is related to the cohomology of the
constant sheaf via appropriate simplicial resolution, which we now
describe.

Let $I_x$ be the indecomposable injective corresponding to $x$. For
a simplex $\sigma$, set $I_\sigma = I_{\min \sigma}$ where $\min
\sigma$ is the minimal element of $\sigma$. If $\tau \subseteq
\sigma$, then $\min \tau \geq \min \sigma$, hence $\Hom_X(I_\tau,
I_\sigma) \simeq k$.

Let $X^{(p)}$ denote the set of $p$-simplices of $X$ and let
$\cI^p_X = \oplus_{\sigma \in X^{(p)}} I_\sigma$. For a $p$-simplex
$\sigma = x_0 < x_1 < \dots < x_p$ and $0 \leq j \leq p$, denote by
$\widehat{\sigma}^j$ the $(p-1)$-simplex obtained from $\sigma$ by
deleting the vertex $x_j$. By considering, for all $\sigma \in
X^{(p)}$ and $0 \leq j \leq p$, the map $I_{\widehat{\sigma}^j} \to
I_\sigma$ corresponding to $(-1)^j \in k$ , we get a map $d^{p-1} :
\cI^{p-1}_X \to \cI^p_X$. The usual sign considerations give $d^p
d^{p-1} = 0$.

\begin{lemma}
\label{l:injsimp} $\hh^i(X) = \hh^i(\Hom_X(k_X,\cI^{\bullet}_X))$
for all $i \geq 0$.
\end{lemma}
\begin{proof}
Indeed, the $p$-th term is $\Hom_X(k_X, \cI^p_X) = \oplus_{\sigma
\in X^{(p)}} \Hom_X(k_X, I_\sigma) \cong \oplus_{\sigma \in X^{(p)}}
k_X(\min \sigma)$ and can be viewed as the space of functions from
$X^{(p)}$ to $k$. Moreover, the differential is exactly the one used
in the definition of simplicial cohomology.
\end{proof}

\begin{lemma}
\label{l:injres} The complex $0 \to k_X \to \cI^0_X
\xrightarrow{d^0} \cI^1_X \xrightarrow{d^1} \dots$ is an injective
resolution of the constant sheaf $k_X$.
\end{lemma}
\begin{proof}
It is enough to check acyclicity at the stalks.

Let $x \in X$. Then $I_\sigma(x) \neq 0$ only if $\min \sigma \geq
x$, hence it is enough to consider the $p$-simplices of $U_x$, and
the complex of stalks at $x$ equals
\[
0 \to k \to \Hom_{U_x}(k_{U_x}, \cI^0_{U_x}) \to \Hom_{U_x}(k_{U_x},
\cI^1_{U_x}) \to \dots
\]
The acyclicity of this complex follows by Lemma~\ref{l:injsimp} with
$X = U_x$, using the fact that $U_x$ has $x$ as the unique minimal
element, hence $\cK(U_x)$ is contractible and $\hh^i(U_x) = 0$ for
$i>0$, $\hh^0(U_x) = k$.
\end{proof}

\begin{prop}
\label{p:simpsheaf} $\hh^i(X;k_X) = \hh^i(X)$ for all $i \geq 0$.
\end{prop}
\begin{proof}
Using Lemma~\ref{l:injsimp} and the injective resolution of
Lemma~\ref{l:injres},
\[
\hh^i(X;k_X) = \hh^i(\Hom_X(k_X, \cI^{\bullet}_X)) = \hh^i(X)
\]
\end{proof}

\subsubsection{Hochschild cohomology}
A $k$-algebra $\gL$ has a natural structure of a
$\gL$-$\gL$-bimodule, or a $\gL \otimes_k \gL^{op}$ right module.
The group $\Ext^i_{\gL \otimes \gL^{op}}(\gL, \gL)$ is called the
\emph{$i$-th Hochschild cohomology} of $\gL$, and we denote it by
$\HH^i(\gL)$.

The Hochschild cohomology of incidence algebras of posets was widely
studied, see~\cite{Cibils89,GaticaRedondo01,GerstenhaberSchack83}.
The following theorem relates the Hochschild cohomology of an
incidence algebra of a poset $X$ with its simplicial cohomology.

\begin{theorem}[\cite{Cibils89,GerstenhaberSchack83}]
\label{t:HC} $\HH^i(kX) = \hh^i(X)$ for all $i \geq 0$.
\end{theorem}

Combining this with Proposition~\ref{p:simpsheaf}, we get:
\begin{cor}
$\HH^i(kX) = \hh^i(X;k_X) = \Ext^i_X(k_X, k_X)$ for all $i \geq 0$.
\end{cor}

\subsubsection{Derived invariants}
\begin{cor}
If $X \sim Y$ then $\beta^i(X) = \beta^i(Y)$ for all $i \geq 0$.
\end{cor}
\begin{proof}
Follows from Theorem~\ref{t:HC} and the fact that the Hochschild
cohomology of a $k$-algebra is preserved under derived
equivalence~\cite{Happel89,Rickard91}.
\end{proof}

The alternating sum $\chi(X) = \sum_{i \geq 0} (-1)^i \beta^i(X)$ is
known as the \emph{Euler characteristic} of $X$.

\begin{cor}
If $X \sim Y$ then $\chi(X) = \chi(Y)$.
\end{cor}

We give two interpretations of $\chi(X)$. First, by
Proposition~\ref{p:simpsheaf},
\[
\chi(X) = \sum_{i \geq 0} (-1)^i \beta^i(X) = \sum_{i \geq 0} \dim_k
\Hom_{\dX}(k_X, k_X[i]) = \bform{[k_X]}{[k_X]}_X
\]
where $[k_X]$ is the image of $k_X$ in $K_0(X)$. Since $[k_X] =
\sum_{x \in X} [S_x]$,
\[
\bform{[k_X]}{[k_X]} = \sum_{x,y \in X} \bform{[S_x]}{[S_y]}_X =
\sum_{x,y \in X} \mu_X(x,y)
\]
hence $\chi(X)$ is the sum of entries of the matrix $\one_X^{-1}$.
We see that not only the $\bZ$-congruence class of $\one_X^{-1}$ is
preserved by derived equivalence, but also the sum of its entries.

For the second interpretation, changing the order of summation we
get
\[
\sum_{x,y \in X} \bform{[S_x]}{[S_y]}_X = \sum_{i \geq 0} (-1)^i
\sum_{x,y \in X} \dim \Ext^i_X(S_x, S_y)
\]
Using the fact that $\dim \Ext^i(S_x, S_y)$ equals $\delta_{xy}$ for
$i=0$; counts the number of arrows from $x$ to $y$ in the Hasse
diagram of $X$ when $i=1$; and counts the number of commutativity
relations between $x$ and $y$ for $i=2$, we see that at least when
$\gldim X \leq 2$, $\chi(X)$ equals the number of points minus the
number of arrows in the Hasse diagram plus the number of relations
etc.

\subsection{Operations preserving derived equivalence}
We show that derived equivalence is preserved under taking opposites
and products.

\begin{defn}
The \emph{opposite} of a poset $X$, denoted by $X^{op}$, is the
poset $(X,\leq^{op})$ with $x \leq^{op} x'$ if and only if $x \geq
x'$.
\end{defn}

\begin{lemma}
Let $\cA$ be an abelian category. Then $Sh_{X^{op}} \cA \simeq (Sh_X
\cA^{op})^{op}$.
\end{lemma}
\begin{proof}
A sheaf $\cF$ over $X^{op}$ with values in $\cA$ is defined via
compatible $\cA$-morphisms between the stalks $\cF(y) \to \cF(x)$
for $x \leq y$. Viewing these morphisms as $\cA^{op}$-morphisms we
identify $\cF$ with a sheaf over $X$ with values in $\cA^{op}$.
Since a morphism of sheaves $\cF \to \cG$ is specified via
compatible $\cA$-morphisms $\cF(x) \to \cG(x)$, this identification
gives an equivalence $Sh_{X^{op}} \cA \simeq (Sh_X \cA^{op})^{op}$.
\end{proof}

\begin{cor}
\label{c:shXop} $Sh_{X^{op}}$ is equivalent to $(Sh_X)^{op}$.
\end{cor}
\begin{proof}
Let $\cA$ be the category of finite dimensional $k$-vector spaces.
Then the functor $V \mapsto V^{\vee}$ mapping a finite dimensional
$k$-vector space to its dual induces an equivalence $\cA \simeq
\cA^{op}$.
\end{proof}

\begin{prop}
If $X \sim Y$ then $X^{op} \sim Y^{op}$.
\end{prop}
\begin{proof}
It is well known that for an abelian category $\cA$, the opposite
category $\cA^{op}$ is also abelian and $\cD^b(\cA) \simeq
\cD^b(\cA^{op})^{op}$ by mapping a complex $K = (K^i)_{i \in \bZ}$
over $\cA$ to the complex $K^{\vee}$ over $\cA^{op}$ with
$(K^{\vee})^i = K^{-i}$.

Applying this for $\cA = Sh_X$ and using Corollary~\ref{c:shXop}, we
deduce that $\cD^b(X^{op}) \simeq \dX^{op}$.
\end{proof}

\begin{defn}
The \emph{product} of two posets $X$, $Y$, denoted $X \times Y$, is
the poset whose underlying set is $X \times Y$, with $(x,y) \leq
(x',y')$ if $x \leq x'$ and $y \leq y'$.
\end{defn}

\begin{lemma}
$k(X \times Y) = kX \otimes_k kY$.
\end{lemma}
\begin{proof}
Observe that the function $kX \otimes_k kY \to k(X \times Y)$
defined by mapping the basis elements $e_{xx'} \otimes e_{yy'}$ to
$e_{(x,y)(x',y')}$ where $x \leq x'$ and $y \leq y'$, is an
isomorphism of $k$-algebras.
\end{proof}

\begin{prop}
If $X_1 \sim X_2$ and $Y_1 \sim Y_2$ then $X_1 \times Y_1 \sim X_2
\times Y_2$.
\end{prop}
\begin{proof}
The claim follows from the previous lemma and the corresponding fact
for tensor products of finite dimensional algebras over $k$,
see~\cite[Lemma~4.3]{Rickard91}.
\end{proof}

\section{Derived equivalences via exceptional collections}

\subsection{Strongly exceptional collections}
Let $k$ be a field and let $\cT$ be a triangulated $k$-category.
\begin{defn}
A sequence $E_1,\dots,E_n$ of objects of $\cT$ is called a
\emph{strongly exceptional collection} if
\begin{align}
\notag
&\Hom_{\cT}(E_s, E_t[i]) = 0 & & 1 \leq s, t \leq n \,,\,i \neq 0 \\
\label{e:SE}
&\Hom_{\cT}(E_s, E_t) = 0 & & 1 \leq s < t \leq n\\
\notag &\Hom_{\cT}(E_s, E_s) = k & & 1 \leq s \leq n
\end{align}

An unordered finite collection $\cE$ of objects of $\cT$ will be
called \emph{strongly exceptional} if it can be ordered in a
sequence which forms a strongly exceptional collection.
\end{defn}

Let $\cE = E_1,\dots,E_n$ be a strongly exceptional collection in
$\cT$, and consider $E = \oplus_{s=1}^{n} E_s$. The
conditions~\eqref{e:SE} imply that $\Hom_{\cT}(E, E[i]) = 0$ for $i
\neq 0$ and that $\End_{\cT}(E)$ is a finite dimensional
$k$-algebra. If $\cE$ generates $\cT$, then $E$ is a \emph{tilting
object} in $\cT$.

For an algebra $A$ over $k$, denote by $\dA$ the bounded derived
category of complexes of finite dimensional right modules over $A$.
The following result of Bondal shows that the existence of a
generating strongly exceptional collection in a derived category
leads to derived equivalence with $\dA$ where $A$ is the
endomorphism algebra of the corresponding tilting object.

\begin{theorem}[\protect{\cite[\S 6]{Bondal90}}]
\label{t:SE} Let $\cA$ be an abelian category and let
$E_1,\dots,E_n$ be a strongly exceptional collection which generates
$\cD^b(\cA)$. Set $E = \oplus_{s=1}^{n} E_s$. Then the functor
\[
\mathbf{R}\mathrm{Hom}(E,-) : \cD^b(\cA) \to \cD^b(\End_{\cD^b(\cA)}
E)
\]
is a triangulated equivalence.
\end{theorem}

When $\cA$ is a category of finite dimensional modules over a finite
dimensional algebra, as in the case of $Sh_X$, the result of the
theorem can also be deduced from Rickard's Morita theory of derived
equivalences of algebras \cite{Rickard89} (see also
\cite[(3.2)]{Keller98}) by observing that $E$ is a so-called
one-sided tilting complex.

\begin{example}
For a poset $X$, the collection $\{P_x\}_{x \in X}$ (and $\{I_x\}_{x
\in X}$) of indecomposable projectives (injectives) is strongly
exceptional, generates $\dX$, and the corresponding endomorphism
algebra is isomorphic to the incidence algebra of $X$ (Use
Lemma~\ref{l:PxPy}).
\end{example}

\subsection{A gluing construction}
Let $\cT, \cT', \cT''$ be three triangulated categories with
triangulated functors
\[
\xymatrix{\cT' \ar@<1ex>[r]^{i_{*}} & \cT \ar@<1ex>[l]^{i^{-1}}
\ar@<1ex>[r]^{j^{-1}} & \cT'' \ar@<1ex>[l]^{j_{!}}}
\]
Assume that there are adjunctions
\begin{align}
\label{e:iadj}
\Hom_{\cT'}(i^{-1}\cF, \cF') &\simeq \Hom_{\cT}(\cF, i_{*}\cF') \\
\label{e:jadj} \Hom_{\cT''}(\cF'', j^{-1}\cF) &\simeq
\Hom_{\cT}(j_{!}\cF'', \cF)
\end{align}
for $\cF \in \cT$, $\cF' \in \cT'$, $\cF'' \in \cT''$. Assume also
that $j^{-1}i_{*} = 0$, $i^{-1}j_{!} = 0$, $i^{-1}i_{*} \simeq
\Id_{\cT'}$ and $j^{-1}j_{!} \simeq \Id_{\cT''}$.

\begin{lemma}
Let $\cF, \cG \in \cT$. Then
\begin{align}
\label{e:iiF} \Hom_{\cT}(\ii\cF, \ii\cG) &\simeq \Hom_{\cT}(\cF,
i_{*}i^{-1}\cG) \\
\label{e:jjF} \Hom_{\cT}(\jj\cF, \jj\cG) &\simeq
\Hom_{\cT}(j_{!}j^{-1}\cF, \cG) \\
\label{e:jiF} \Hom_{\cT}(\jj\cF, \ii\cG) &= 0
\end{align}
\end{lemma}
\begin{proof}
The claims follow from the adjunctions~\eqref{e:iadj},\eqref{e:jadj}
and our additional hypotheses. For example, for the first claim
use~\eqref{e:iadj} and $i^{-1}i_{*} \simeq \Id_{\cT'}$ to get that
\[
\begin{split}
\Hom_{\cT}(\ii\cF, \ii\cG) \simeq \Hom_{\cT'}(i^{-1} \ii\cF,
i^{-1}\cG) = \\
= \Hom_{\cT'}(i^{-1}\cF, i^{-1}\cG) \simeq \Hom_{\cT}(\cF, \ii\cG)
\end{split}
\]
\end{proof}

We apply this for the following situation, cf.~\cite[\S1.4]{BBD82}.
Let $X$ be a poset and let $Y \subseteq X$ be a closed subset, $U =
X \setminus Y$ its complement. Denote by $i: Y \to X$, $j:U \to X$
the embeddings. Since the functors $i_*, j_!$ are exact, we can
consider the functors
\begin{align*}
i^{-1} : \dX \to \dY & & j^{-1} : \dX \to \dU \\
i_{*} : \dY \to \dX & & j_{!} : \dU \to \dX
\end{align*}
between the derived categories. Taking $\cT = \dX$, $\cT' = \dY$ and
$\cT'' = \dU$, we see that the above assumptions are satisfied,
where the adjunctions~\eqref{e:iadj}, \eqref{e:jadj} follow
from~\eqref{e:adjder}.

For $y \in Y$ and $u \in U$, let $\wtP_y = \ii P_y$ and $\wtI_u =
\jj I_u$ be ``truncated'' versions of the projectives and
injectives. Explicitly,
\[
\wtP_y(x) = \begin{cases} k & x \in Y \,,\, y \leq x \\ 0 &
\text{otherwise} \end{cases} \qquad \wtI_u(x) = \begin{cases} k & x
\in U \,,\, x \leq u \\ 0 & \text{otherwise} \end{cases}
\]
with identity maps between nonzero stalks.

\begin{prop}
\label{p:PISE} The collection $\cE_Y = \left\{ \wtP_y \right\}_{y
\in Y} \bigcup \left\{ \wtI_u[1] \right\}_{u \in U} $ is strongly
exceptional and generates $\dX$.
\end{prop}
\begin{proof}
Let $y, y' \in Y$. By~\eqref{e:iiF},
\begin{equation}
\label{e:homPP} \Hom(\wtP_y, \wtP_{y'}) \simeq \Hom(P_y, \wtP_{y'})
= \wtP_{y'}(y) =
\begin{cases} k &\text{if $y' \leq y$} \\ 0 & \text{otherwise}
\end{cases}
\end{equation}
and $\Hom(\wtP_y, \wtP_{y'}[n]) = 0$ for $n \neq 0$. Similarly, for
$u, u' \in U$, by~\eqref{e:jjF},
\begin{equation}
\label{e:homII} \Hom(\wtI_u, \wtI_{u'}) \simeq \Hom(\wtI_u, I_{u'})
= \wtI_u (u') =
\begin{cases} k &\text{if $u' \leq u$} \\ 0 & \text{otherwise}
\end{cases}
\end{equation}
and $\Hom(\wtI_u, \wtI_{u'}[n]) = 0$ for $n \neq 0$.

Let $y \in Y$ and $u \in U$. By~\eqref{e:jiF}, $\Hom(\wtI_u,
\wtP_y[n]) = 0$ for all $n \in \bZ$. Consider now $\Hom(\wtP_y,
\wtI_u[n])$. The distinguished triangle $\wtI_u \to I_u \to \ii I_u
\to \wtI_u[1]$ of~\eqref{e:jjii} gives rise to a long exact sequence
\begin{equation}
\label{e:PyIuLES} \dots \to \Hom(\wtP_y, \wtI_u) \to \Hom(\wtP_y,
I_u) \to \Hom(\wtP_y, \ii I_u) \to \dots
\end{equation}
Since $I_u$ is injective, $\Hom(\wtP_y, I_u[n]) = 0$ for $n \neq 0$
and $\Hom(\wtP_y, I_u) = \wtP_y(u) = 0$. Therefore \eqref{e:PyIuLES}
induces isomorphisms
\begin{equation}
\label{e:PIi} \Hom(\wtP_y, \ii I_u[n]) \xrightarrow{\simeq}
\Hom(\wtP_y, \wtI_u[n+1])
\end{equation}
for all $n \in \bZ$. By~\eqref{e:iiF},
\[
\Hom(\wtP_y, \ii I_u[n]) = \Hom(P_y, \ii I_u[n]) = \begin{cases}
(\ii I_u)(y) & n = 0 \\ 0 & n \neq 0
\end{cases}
\]
and $(\ii I_u)(y) = k$ if $y < u$ and $0$ otherwise, hence
\begin{equation}
\label{e:homPI} \Hom(\wtP_y, \wtI_u[1]) = \begin{cases} k & \text{if $y < u$} \\
0 & \text{otherwise} \end{cases}
\end{equation}
and $\Hom(\wtP_y, \wtI_u[1+n]) = 0$ for $n \neq 0$.

Note that one can also compute $\Hom(\wtP_y, \wtI_u[n])$ by
considering the triangle $\jj P_y \to P_y \to \wtP_y \to \jj P_y[1]$
and using the induced isomorphisms
\begin{equation}
\label{e:jPI} \Hom(\jj P_y, \wtI_u[n]) \xrightarrow{\simeq}
\Hom(\wtP_y, \wtI_u[n+1])
\end{equation}

The above calculations show that if we order each of the sets $Y$
and $U$ linearly extending the partial order induced by $X$ and
arrange the elements of $\cE_Y$ in a sequence by first taking the
elements of $U$ and then taking those of $Y$, we get a strongly
exceptional collection.

To prove that $\cE_Y$ generates $\dX$, it is enough to show that
every sheaf belongs to the triangulated subcategory generated by
$\cE_Y$. By~\eqref{e:jjii}, it is enough to verify this for
$i_{*}\cF'$ and $j_{!} \cF''$ where $\cF' \in Sh_Y$, $\cF'' \in
Sh_U$. The collection of sheaves $i^{-1}P_y$, being a complete set
of indecomposable projectives of $Sh_Y$, generates $\dY$. Similarly,
the sheaves $j^{-1}I_u$ form a complete set of indecomposable
injectives of $Sh_U$ and generate $\dU$. Now the result follows by
applying the triangulated functors $i_*, j_!$.
\end{proof}

\subsection{The endomorphism algebras $A_Y$}
Fix a poset $X$, and let $Y \subseteq X$ be a closed subset.
Consider $T_Y = (\oplus_{y \in Y} \wtP_y) \oplus (\oplus_{u \in U}
\wtI_u)[1]$ and let $A_Y = \End_{\dX} T_Y$. By Theorem~\ref{t:SE}
and Proposition~\ref{p:PISE}, we have:
\begin{cor}
\label{c:dXAY} $\dX \simeq \cD^b(A_Y)$.
\end{cor}

\begin{prop}
The algebra $A_Y$ has as a $k$-basis the elements
\[
\left\{ e_{yy'} \,:\, y \leq y' \right\} \cup \left\{e_{u'u} \,:\,
u' \leq u \right\} \cup \left\{ e_{uy} \,:\, y < u \right\}
\]
where $y,y' \in Y$, $u',u \in U$. The multiplication is defined by
\begin{align*}
&e_{y y'} e_{y' y''} = e_{yy''} & &
e_{u''u'} e_{u'u} = e_{u''u} \\
&e_{uy} e_{y y'} = \begin{cases} e_{u y'} & y' < u \\ 0 &
\text{otherwise} \end{cases} & & e_{u'u} e_{uy} = \begin{cases}
e_{u'y} & y < u' \\ 0 & \text{otherwise} \end{cases}
\end{align*}
for $y \leq y' \leq y'' \in Y$, $u'' \leq u' \leq u \in U$ (all
other products are zero).
\end{prop}
\begin{proof}
For $y \leq y' \in Y$, using~\eqref{e:homPP}, choose $e_{yy'} \in
\Hom(\wtP_{y'}, \wtP_y)$ corresponding to $1 \in \wtP_y(y')$. In
other words, the stalk of the morphism $e_{yy'}$ at $y'$ is the
identity map on $k$. Then $e_{yy'}e_{y'y''}=e_{yy''}$ for $y \leq y'
\leq y'' \in Y$.

Similarly, for $u' \leq u \in U$, using~\eqref{e:homII}, choose
$e_{u'u} \in \Hom(\wtI_u[1], \wtI_{u'}[1])$ corresponding to $1 \in
\wtI_u(u')$. The stalk of $e_{u'u}$ at $u'$ is the identity map on
$k$ and we have $e_{u''u'}e_{u'u}=e_{u''u}$ for all $u'' \leq u'
\leq u \in U$.

Now consider $y \in Y$ and $u \in U$ such that $y < u$. Using the
isomorphisms~\eqref{e:PIi} and~\eqref{e:jPI}, we have
\[
\xymatrix{ \Hom(\jj P_y, \wtI_u) \ar[d]_{\simeq} \ar[r]^{\simeq} &
\Hom(\wtP_y, \wtI_u[1]) & \ar[l]_{\simeq} \Hom(\wtP_y, \ii I_u)
\ar[d]^{\simeq} \\ (\jj P_y)(u) & & (\ii I_u)(y) }
\]

There are unique $e_{uy}, \widetilde{e}_{uy} \in \Hom(\wtP_y,
\wtI_u[1])$ such that the image of $e_{uy}$ in $(\jj P_y)(u)$ equals
$1$ and the image of $\widetilde{e}_{uy}$ in $(\ii I_u)(y)$ is $1$.
The formula for $e_{u'u}e_{uy}$ where $u' \leq u$ now follows by
considering the composition
\[
\xymatrix{ \Hom(\jj P_y, \wtI_u) \ar[d]_{e_{u'u} \circ -}
\ar[r]^{\simeq} & \Hom(\wtP_y, \wtI_u[1]) \ni e_{uy} \\
\Hom(\jj P_y, \wtI_{u'}) \ar[r]^{\simeq} & \Hom(\wtP_y,
\wtI_{u'}[1]) \ni e_{u'y}}
\]

The formula for $e_{uy} e_{yy'}$ would follow in a similar manner by
considering the composition $- \circ e_{yy'}$ once we know that the
scalar ratio between $\widetilde{e}_{uy}$ and $e_{uy}$ is
\emph{independent} of $u$ and $y$.

Indeed, replacing the objects $\wtP_y$ and $\wtI_u[1]$ by the
quasi-isomorphic complexes $(\jj P_y \to P_y)[1]$ and $(I_u \to \ii
I_u)[1]$, we see that $\Hom(\wtP_y, \wtI_u[1])$ equals the set of
morphisms $(\lambda, \mu)$ between the two complexes
\[
\xymatrix{ {\dots} \ar[r] & 0 \ar[r]  & \jj P_y \ar[r]
\ar[d]^{\lambda} & P_y \ar[r] \ar[d]^{\mu} \ar@{--}[dl] & 0 \ar[r]  & {\dots} \\
{\dots} \ar[r] & 0 \ar[r] & I_u \ar[r] & \ii I_u \ar[r] & 0 \ar[r] &
{\dots} }
\]
modulo homotopy. Note that $(\lambda,\mu) \sim (\lambda',\mu')$ if
and only if $\lambda - \mu = \lambda'-\mu'$. The morphism $e_{uy}$
corresponds to the pair $(1,0)$ while $\widetilde{e}_{uy}$
corresponds to $(0,1)$, hence $\widetilde{e}_{uy} = -e_{uy}$.

It is clear that the elements constructed above form a $k$-basis of
$A_Y$ and satisfy the required relations.
\end{proof}

\begin{example}
Let $X$ be the poset with Hasse diagram as in the left picture, and
let $Y = \{1\}$. The algebra $A_Y$ is shown in the right picture, as
the path algebra of the quiver $A_3$ modulo the zero relation
indicated by the dotted arrow (i.e. the product of $2 \to 3$ and $3
\to 1$ is zero).
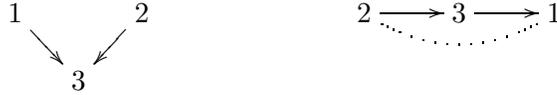
\begin{align*}
\xymatrix@=1pc{1 \ar[dr] & & 2 \ar[dl] \\ & 3} & & \xymatrix{2
\ar@/_1pc/@{.}[rr] \ar[r] & 3 \ar[r] & 1}
\end{align*}
\end{example}

\begin{lemma}
\label{l:AYposet} Let $X' = U \cup Y$ and define a binary relation
$\leq'$ on $X'$ by
\begin{equation}
\label{e:X'poset} u' \leq' u \Leftrightarrow u' \leq u \quad y \leq'
y' \Leftrightarrow y \leq y' \quad u <' y \Leftrightarrow y < u
\end{equation}
for $u,u' \in U$, $y,y' \in Y$. Then $\leq'$ is a partial order if
and only if the following condition holds:
\begin{equation}
\tag{$\star$} \label{e:AYposet} \text{Whenever $y \leq y' \in Y$,
$u' \leq u \in U$ and $y < u$, we have that $y' < u'$.}
\end{equation}

When this condition holds, the endomorphism algebra $A_Y$ is
isomorphic to the incidence algebra of $(X',\leq')$.
\end{lemma}
\begin{proof}
The first part is clear from the requirement of transitivity of
$\leq'$.

The condition~\eqref{e:AYposet} implies that $e_{u'u}e_{uy} =
e_{u'y}$ and $e_{uy}e_{yy'} = e_{uy'}$ whenever $u' \leq u$, $y \leq
y'$ and $y < u$, so that $A_Y$ is the incidence algebra of
$(X',\leq')$.
\end{proof}

\subsection{Lexicographic sums along bipartite graphs}

\begin{defn}
Let $S$ be a poset, and let $\fX = \{X_s\}_{s \in S}$ be a
collection of posets indexed by the elements of $S$. The
\emph{lexicographic sum of $\fX$ along $S$}, denoted $\oplus_S \fX$,
is the poset $(X,\leq)$ where $X = \coprod_{s \in S} X_s$ is the
disjoint union of the $X_s$ and for $x \in X_s$, $y \in X_t$ we have
$x \leq y$ if either $s<t$ (in $S$) or $s=t$ and $x \leq y$ (in
$X_s$).
\end{defn}

\begin{example} 
The usual \emph{ordinal sum} $X_1 \oplus X_2 \oplus \dots \oplus
X_n$ of $n$ posets is the lexicographic sum of $\{X_1,\dots,X_n\}$
along the chain $1 < 2 < \dots < n$.
\end{example}

\begin{defn}
A poset $S$ is called a \emph{bipartite graph} if it can be written
as a disjoint union of two nonempty subsets $S_0$ and $S_1$ such
that $s < s'$ in $S$ implies that $s \in S_0$ and $s' \in S_1$.
\end{defn}

It follows from the definition that the posets $S_0$, $S_1$ are
\emph{anti-chains}, that is, no two distinct elements in $S_0$ (or
$S_1$) are comparable.

\begin{example}
The left Hasse diagram represents a bipartite poset $S$. The right
one is the Hasse diagram of its opposite $S^{op}$.
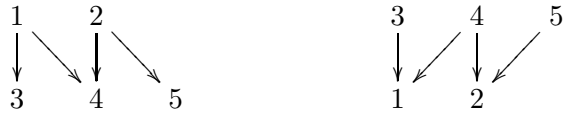
\begin{align*}
\xymatrix@=1.5pc{
1 \ar[d] \ar[dr] & 2 \ar[d] \ar[dr] \\
3 & 4 & 5} & &
\xymatrix@=1.5pc{
3 \ar[d] & 4 \ar[dl] \ar[d] & 5 \ar[dl] \\
1 & 2 }
\end{align*}

Let $\fX = \left\{ X_1, X_2, X_3, X_4, X_5 \right\}$ be the
collection
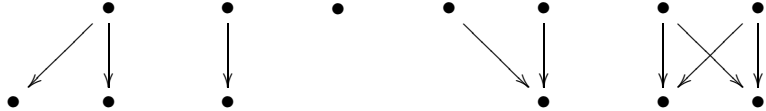
\begin{align*}
\xymatrix{ & {\bullet} \ar[dl] \ar[d] \\ {\bullet} & {\bullet}} &
& \xymatrix{ {\bullet} \ar[d] \\ {\bullet}} &
& {\bullet} &
& \xymatrix{ {\bullet} \ar[dr] & {\bullet} \ar[d] \\ & {\bullet}} &
& \xymatrix{ {\bullet} \ar[dr] \ar[d] & {\bullet} \ar[dl] \ar[d] \\
{\bullet} & {\bullet}}
\end{align*}

The graphs shown below are the Hasse diagrams of $\oplus_S \fX$
(left) and $\oplus_{S^{op}} \fX$ (right).
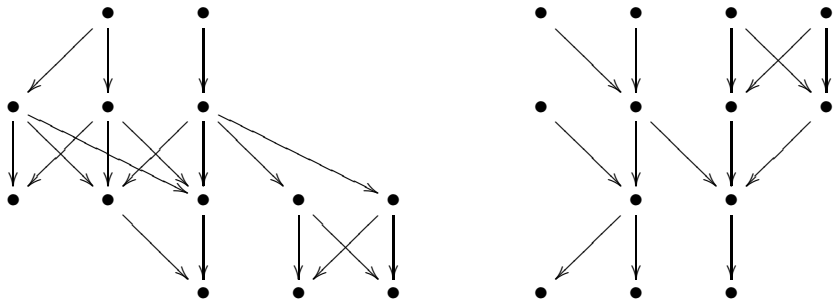
\[
\xymatrix{
& \bullet \ar[dl] \ar[d] & \bullet \ar[d] \\
\bullet \ar[d] \ar[dr] \ar[drr]& \bullet \ar[dl] \ar[d] \ar[dr] &
\bullet
\ar[dl] \ar[d] \ar[dr] \ar[drr] \\
\bullet & \bullet \ar[dr] & \bullet \ar[d] & \bullet \ar[d]
\ar[dr] & \bullet \ar[dl] \ar[d] \\
& & \bullet & \bullet & \bullet} \qquad \qquad
\xymatrix{
 \bullet \ar[dr] & \bullet \ar[d] & \bullet \ar[d] \ar[dr] & \bullet \ar[dl] \ar[d] \\
 \bullet \ar[dr] & \bullet \ar[d] \ar[dr] & \bullet \ar[d] & \bullet \ar[dl] \\
& \bullet \ar[dl] \ar[d] & \bullet \ar[d] \\
\bullet & \bullet & \bullet}
\]
\end{example}

\begin{theorem}
\label{t:bipartite} If $S$ is a bipartite graph and
$\fX=\{X_s\}_{s\in S}$ is a collection of posets, then $\oplus_S \fX
\sim \oplus_{S^{op}} \fX$.
\end{theorem}
\begin{proof}
Let $S = S_0 \amalg S_1$ be a partition as in the definition of
bipartite poset. Let $\fX_0 = \{X_s\}_{s \in S_0}$, $\fX_1 =
\{X_s\}_{s \in S_1}$ and let $X = \oplus_S \fX$, $Y = \oplus_{S_0}
\fX_0$, $U = \oplus_{S_1} \fX_1$. The sets $Y$ and $U$ can be viewed
as disjoint subsets of $X$ with $X = Y \cup U$. Moreover, since
there are no relations $s_1 < s_0$ with $s_0 \in S_0$, $s_1 \in
S_1$, there are no relations $u < y$ with $y \in Y$, $u \in U$, thus
$Y$ is closed and $U$ is open in $X$. By Corollary~\ref{c:dXAY},
$\dX \simeq \cD^b(A_Y)$ where $A_Y$ is the endomorphism algebra of
the direct sum of the strongly exceptional collection of
Proposition~\ref{p:PISE}.

We show that the condition~\eqref{e:AYposet} of
Lemma~\ref{l:AYposet} holds. Indeed, let $y \leq y' \in Y$, $u' \leq
u \in U$. There exist $s_0, s_0' \in S_0$, $s_1, s_1' \in S_1$ such
that $y \in X_{s_0}$, $y' \in X_{s_0'}$, $u \in X_{s_1}$ and $u' \in
X_{s'_1}$. Now, $s_0' = s_0$ and $s_1' = s_1$ since $y \leq y'$, $u'
\leq u$ and $S_0$, $S_1$ are anti-chains. If $y < u$, then $s_0 <
s_1$, hence $y' < u'$ and~\eqref{e:AYposet} is satisfied. Therefore
$A_Y$ is the incidence algebra of the poset $X'$ defined
in~\eqref{e:X'poset}.

Since $X'$ is a disjoint union of the posets $U$ and $Y$ with the
original order inside each but with reverse order between them, it
is easy to see that $X'$ equals the lexicographic sum of $\fX$ along
the opposite poset $S^{op}$.
\end{proof}

\begin{cor}
\label{c:XplusY} Let $X, Y$ be two posets. Then $X \oplus Y \sim Y
\oplus X$
\end{cor}
\begin{proof}
Take $S$ to be the chain $1 < 2$.
\end{proof}

As special cases, we obtain the following two well known examples.

\begin{example}
The following two posets (represented by their Hasse diagrams) are
derived equivalent.
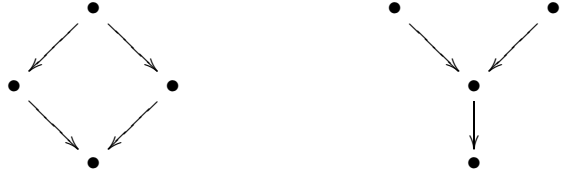
\begin{align*}
\xymatrix@=1.5pc{
& {\bullet} \ar[dl] \ar[dr] \\
{\bullet} \ar[dr] & & {\bullet} \ar[dl] \\
& {\bullet} } & &
\xymatrix@=1.5pc{
{\bullet} \ar[dr] & & {\bullet} \ar[dl] \\
& {\bullet} \ar[d] \\
& {\bullet} }
\end{align*}
The right poset is obtained from the left one by an APR
tilt~\cite{APR79}, see also \cite[(III.2.14)]{Happel88}.
\end{example}

\begin{example}
The two posets below are derived equivalent.
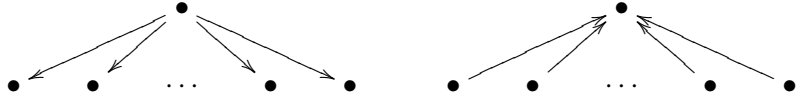
\begin{align*}
\xymatrix@=1.5pc{ & & {\bullet} \ar[dll] \ar[dl] \ar[dr] \ar[drr] \\
{\bullet} & {\bullet} & {\cdots} & {\bullet} & {\bullet}} & &
\xymatrix@=1.5pc{ & & {\bullet} \\
{\bullet} \ar[urr] & {\bullet} \ar[ur] & {\cdots} & {\bullet}
\ar[ul] & {\bullet} \ar[ull]}
\end{align*}
This is a special case of BGP reflection \cite{BGP73}, turning a
source into a sink (and vice versa).
\end{example}

\begin{cor}
Let $S$ be a bipartite graph. Then $S \sim S^{op}$.
\end{cor}
\begin{proof}
Take in Theorem~\ref{t:bipartite} each $X_s$ to be a point.
\end{proof}

Note that the last Corollary can also be deduced from \cite{BGP73}
since $Sh_S$ is the category of representations of a quiver without
oriented cycles, namely the Hasse diagram of $S$, and $S^{op}$ is
obtained from $S$ by reverting all the arrows.

\subsection{Ordinal sums of three posets}
The result of Corollary~\ref{c:XplusY} raises the natural question
whether the derived equivalence class of an ordinal sum of more than
two posets does not depend on the order of the summands. The
following proposition shows that it is enough to consider the case
of three summands.

\begin{prop}
Let $\cX$ be a family of posets closed to taking ordinal sums.
Assume that for any three posets $X, Y, Z \in \cX$,
\begin{equation}
\label{e:XYZ} X \oplus Y \oplus Z \sim Y \oplus X \oplus Z
\end{equation}

Then for any $n \geq 1$, $\pi \in S_n$ and $X_1,\dots,X_n \in \cX$,
\[
X_{\pi(1)} \oplus \dots \oplus X_{\pi(n)} \sim X_1 \oplus \dots
\oplus X_n
\]
\end{prop}
\begin{proof}
For $n=1$ the claim is trivial and for $n=2$ it is just
Corollary~\ref{c:XplusY}. Let $n \geq 3$ and consider the set $G_n$
of permutations in $\pi \in S_n$ such that $X_{\pi(1)} \oplus \dots
\oplus X_{\pi(n)} \sim X_1 \oplus \dots \oplus X_n$ for all
$X_1,\dots,X_n \in \cX$. Then $G_n$ is a subgroup of $S_n$, and the
claim to be proved is that $G_n = S_n$.

Let $X_1,\dots,X_n \in \cX$. Taking $X=X_1$ and $Y=X_2 \oplus \dots
\oplus X_n$, we see by Corollary~\ref{c:XplusY} that the cycle
$(1\,2\,\dots\,n)$ belongs to $G_n$. Now take $X=X_1$, $Y=X_2$ and
$Z= X_3 \oplus \dots \oplus X_n$. By~\eqref{e:XYZ}, $Y \oplus X
\oplus Z \sim X \oplus Y \oplus Z$, hence $(1\,2) \in G_n$. The
claim now follows since $(1\,2)$ and $(1\,2\,\dots\,n)$
generate~$S_n$.
\end{proof}

We give a counterexample to show that~\eqref{e:XYZ} is false in
general.

\begin{figure}
\begin{align*}
\xymatrix{ \bullet \ar[d] \ar[dr] & \bullet \ar[dl] \ar[d] & \bullet
\ar[dll] \ar[dl] \\
\bullet \ar[d] & \bullet \ar[ddl] \ar[dd] \ar[ddr] \\
\bullet \ar[d] \ar[dr] \ar[drr] \\
\bullet \ar[d] \ar[dr] & \bullet \ar[dl] \ar[d] & \bullet
\ar[dll] \ar[dl] \\
\bullet \ar[d] & \bullet \\
\bullet } & & \xymatrix{
\bullet \ar[d] & \bullet \ar[ddl] \ar[dd] \ar[ddr] \\
\bullet \ar[d] \ar[dr] \ar[drr] \\
\bullet \ar[d] \ar[dr] \ar[drr] & \bullet \ar[dl] \ar[d] \ar[dr] &
\bullet \ar[dll] \ar[dl] \ar[d] \\
\bullet \ar[d] \ar[dr] & \bullet \ar[dl] \ar[d] & \bullet
\ar[dll] \ar[dl] \\
\bullet \ar[d] & \bullet \\
\bullet } \\
X \oplus Y \oplus Z & & Y \oplus X \oplus Z
\end{align*}
\caption{Two posets which are not derived equivalent despite their
structure as ordinal sums of the same three posets in different
orders.} \label{f:posets}
\end{figure}
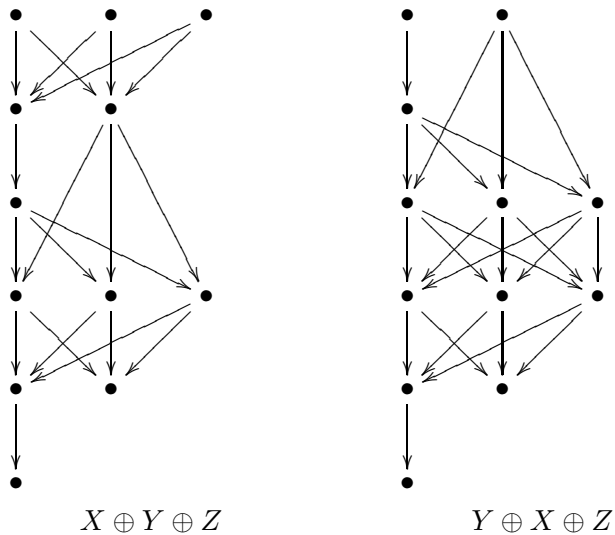

\begin{example}
Let
\begin{align*}
X = \xymatrix@=1.5pc{\bullet & \bullet & \bullet} & &
Y = \xymatrix@=1.5pc{\bullet \ar[d] & \bullet \\ \bullet}
\end{align*}
and let $Z = X \oplus Y$. Then the posets $X \oplus Y \oplus Z$ and
$Y \oplus X \oplus Z$, depicted in Figure~\ref{f:posets}, are not
derived equivalent since their Euler forms are not equivalent over
$\bZ$ (they are equivalent over $\bQ$, though). This is shown using
Corollary~\ref{c:simQFp} with the prime $p=11$.
\end{example}


\begin{thebibliography}{10}

\bibitem{APR79}
{\sc Auslander, M., Platzeck, M.~I., and Reiten, I.}
\newblock Coxeter functors without diagrams.
\newblock {\em Trans. Amer. Math. Soc. 250\/} (1979), 1--46.

\bibitem{BBD82}
{\sc Be{\u\i}linson, A.~A., Bernstein, J., and Deligne, P.}
\newblock Faisceaux pervers.
\newblock In {\em Analysis and topology on singular spaces, I (Luminy, 1981)},
  vol.~100 of {\em Ast\'erisque}. Soc. Math. France, Paris, 1982, pp.~5--171.

\bibitem{BGP73}
{\sc Bern{\v{s}}te{\u\i}n, I.~N., Gel{\cprime}fand, I.~M., and
Ponomarev,
  V.~A.}
\newblock Coxeter functors, and {G}abriel's theorem.
\newblock {\em Uspehi Mat. Nauk 28}, 2(170) (1973), 19--33.

\bibitem{Bondal90}
{\sc Bondal, A.}
\newblock {Representation of associative algebras and coherent sheaves.}
\newblock {\em Math. USSR, Izv. 34}, 1 (1990), 23--42.

\bibitem{BondalOrlov02}
{\sc Bondal, A., and Orlov, D.}
\newblock Derived categories of coherent sheaves.
\newblock In {\em Proceedings of the International Congress of Mathematicians,
  Vol. II (Beijing, 2002)\/} (Beijing, 2002), Higher Ed. Press, pp.~47--56.

\bibitem{Cibils89}
{\sc Cibils, C.}
\newblock Cohomology of incidence algebras and simplicial complexes.
\newblock {\em J. Pure Appl. Algebra 56}, 3 (1989), 221--232.

\bibitem{DGM00}
{\sc Deligne, P., Goresky, M., and MacPherson, R.}
\newblock L'alg\`ebre de cohomologie du compl\'ement, dans un espace affine,
  d'une famille finie de sous-espaces affines.
\newblock {\em Michigan Math. J. 48\/} (2000), 121--136.

\bibitem{DuggerShipley04}
{\sc Dugger, D., and Shipley, B.}
\newblock {$K$}-theory and derived equivalences.
\newblock {\em Duke Math. J. 124}, 3 (2004), 587--617.

\bibitem{GaticaRedondo01}
{\sc Gatica, M.~A., and Redondo, M.~J.}
\newblock Hochschild cohomology and fundamental groups of incidence algebras.
\newblock {\em Comm. Algebra 29}, 5 (2001), 2269--2283.

\bibitem{GerstenhaberSchack83}
{\sc Gerstenhaber, M., and Schack, S.~D.}
\newblock Simplicial cohomology is {H}ochschild cohomology.
\newblock {\em J. Pure Appl. Algebra 30}, 2 (1983), 143--156.

\bibitem{Happel88}
{\sc Happel, D.}
\newblock {\em Triangulated categories in the representation theory of
  finite-dimensional algebras}, vol.~119 of {\em London Mathematical Society
  Lecture Note Series}.
\newblock Cambridge University Press, Cambridge, 1988.

\bibitem{Happel89}
{\sc Happel, D.}
\newblock Hochschild cohomology of finite-dimensional algebras.
\newblock In {\em S\'eminaire d'Alg\`ebre Paul Dubreil et Marie-Paul Malliavin,
  39\`eme Ann\'ee (Paris, 1987/1988)}, vol.~1404 of {\em Lecture Notes in
  Math.} Springer, Berlin, 1989, pp.~108--126.

\bibitem{IgusaZacharia90}
{\sc Igusa, K., and Zacharia, D.}
\newblock On the cohomology of incidence algebras of partially ordered sets.
\newblock {\em Comm. Algebra 18}, 3 (1990), 873--887.

\bibitem{Keller98}
{\sc Keller, B.}
\newblock On the construction of triangle equivalences.
\newblock In {\em Derived equivalences for group rings}, vol.~1685 of {\em
  Lecture Notes in Math.} Springer, Berlin, 1998, pp.~155--176.

\bibitem{McCord66}
{\sc McCord, M.~C.}
\newblock Singular homology groups and homotopy groups of finite topological
  spaces.
\newblock {\em Duke Math. J. 33\/} (1966), 465--474.

\bibitem{Mitchell68}
{\sc Mitchell, B.}
\newblock On the dimension of objects and categories. {II}. {F}inite ordered
  sets.
\newblock {\em J. Algebra 9\/} (1968), 341--368.

\bibitem{Quillen73}
{\sc Quillen, D.}
\newblock Higher algebraic {$K$}-theory. {I}.
\newblock In {\em Algebraic $K$-theory, I: Higher $K$-theories (Proc. Conf.,
  Battelle Memorial Inst., Seattle, Wash., 1972)}. Springer, Berlin, 1973,
  pp.~85--147. Lecture Notes in Math., Vol. 341.

\bibitem{Rickard89}
{\sc Rickard, J.}
\newblock Morita theory for derived categories.
\newblock {\em J. London Math. Soc. (2) 39}, 3 (1989), 436--456.

\bibitem{Rickard91}
{\sc Rickard, J.}
\newblock Derived equivalences as derived functors.
\newblock {\em J. London Math. Soc. (2) 43}, 1 (1991), 37--48.

\bibitem{Spears72}
{\sc Spears, W.~T.}
\newblock Global dimension in categories of diagrams.
\newblock {\em J. Algebra 22\/} (1972), 219--222.

\bibitem{Stong66}
{\sc Stong, R.~E.}
\newblock Finite topological spaces.
\newblock {\em Trans. Amer. Math. Soc. 123\/} (1966), 325--340.

\end{thebibliography}

\def\cprime{$'$}

\end{document}